\def\vers{Feb.~21, 2008, v.2}
\magnification=1200
\hsize=6.5truein
\vsize=8.9truein
\font\bigfont=cmr10 at 14pt
\font\mfont=cmr9

\font\mbfont=cmbx9

\def\scirc{\,\raise.2ex\hbox{${\scriptstyle\circ}$}\,}
\def\mopl{\hbox{$\bigoplus$}}
\def\msum{\hbox{$\sum$}}
\def\C{{\bf C}}
\def\D{{\bf D}}

\def\P{{\bf P}}
\def\Q{{\bf Q}}
\def\R{{\bf R}}
\def\Z{{\bf Z}}
\def\H{{\bf H}}
\def\cH{{\cal H}}
\def\M{{\cal M}}
\def\O{{\cal O}}
\def\cP{{\cal P}}
\def\tC{\widetilde{C}}
\def\tN{\widetilde{N}}
\def\tX{\widetilde{X}}
\def\tXC{\widetilde{X_C}}
\def\tpi{\widetilde{\pi}}
\def\ok{\bar{k}}
\def\End{\hbox{\rm End}}
\def\Ext{\hbox{\rm Ext}}
\def\Hom{\hbox{\rm Hom}}
\def\CH{\hbox{\rm CH}}
\def\Cor{\hbox{\rm Cor}}
\def\Alb{\hbox{\rm Alb}}
\def\Pic{\hbox{\rm Pic}}
\def\NS{\hbox{\rm NS}}
\def\IC{\hbox{\rm IC}}
\def\Prym{\hbox{\rm Prym}}
\def\Sing{\hbox{\rm Sing}\,}
\def\MHM{{\rm MHM}}
\def\MHS{{\rm MHS}}
\def\alg{{\rm alg}}
\def\BM{{\rm BM}}
\def\Gal{{\rm Gal}}
\def\char{{\rm char}\,}
\def\Spec{{\rm Spec}\,}
\def\simto{\buildrel\sim\over\rightarrow}
\def\lra{\longrightarrow}

\hbox{}
\vskip 1cm

\centerline{\bigfont Relative Chow-K\"unneth decompositions}

\smallskip
\centerline{\bigfont for conic bundles and Prym varieties}

\bigskip
\centerline{Jan Nagel and Morihiko Saito}

\bigskip\bigskip
{\narrower\noindent
{\mbfont Abstract.} {\mfont
We construct a relative Chow-K\"unneth decomposition for a conic
bundle over a surface such that the middle projector gives the Prym
variety of the associated double covering of the discriminant of the
conic bundle. This gives a refinement (up to an isogeny) of
Beauville's theorem on the relation between the intermediate Jacobian
of the conic bundle and the Prym variety of the double covering.}
\par}

\bigskip\bigskip
\centerline{\bf Introduction}

\bigskip\noindent
Let $f:X\to S$ be a conic bundle over a surface, i.e.,
$X$ is a smooth projective threefold over $k$, $S$ is a projective
surface over $k$ and the fibers of $f$ are conics, where $k$ is a
perfect field with $\char k\ne 2$.
Let $C$ be the discriminant of $f$; it is a curve whose singularities
are ordinary double points, see [3].
(Here $C$ is not necessarily connected.)
The singularities of $C$ are the points $s\in S$ such that $f^{-1}(s)$
is a double line.
Put $X_C=f^{-1}(C)$, and let $\tXC$ be its normalization (which is
smooth).
Let $\tC$ denote $F_1(X_C/C)$, the relative Fano scheme of lines of
$X_C$ over $C$ (i.e. its fiber over $s\in C$ consists of the
irreducible components of $f^{-1}(s))$.
In [3, 0.3] it is shown that the canonical morphism $\rho:\tC\to C$
is an admissible double covering ('pseudo--rev\^etement' in the
terminology of [3]). Hence $\rho:\tC\to C$ is an \'etale double
covering outside $\Sing(C)$ and the inverse image of a double point
of $C$ is an ordinary double point of $\tC$.
Let $D$ and $C^{\prime}$ denote respectively the normalizations of
$\tC$ and $C$ (which are denoted respectively by $\tN$ and $N$ in
[3], [6]).
Let $C_j$ be the irreducible components of $C$. Let $C_j^{\prime}$ be
the normalization of $C_j$, and $D_j$ be the union of the irreducible
components of $D$ whose image in $C$ is $C_j$.
Renumbering the $C_j$ if necessary, there are integers $r\ge r'\ge 0$
such that the restriction $\rho'_j$ of the double covering
$\rho:\tC\to C$ over $C_j\setminus \Sing C$ is trivial if and only
if $1\le j\le r'$, and the base change of $\rho'_j$ by $k\to\ok$ is
trivial if and only if $1\le j\le r$.

Let $\cP_X$ be the generalized Prym variety associated to the
double covering $\rho:\tC\to C$, as defined in [3, 0.3.2].
Then $\cP_X$ is isogenous to the product of the Prym varieties of
$D_j/C'_j$, see [3], Prop.~0.3.3 (cf. also [6, Prop. 1.5]).
Let $\sigma_j$ be the involution of $D_j$ associated to the double
covering
$$\rho_j: D_j\to C_j^{\prime}.$$
Identifying $\sigma_j$ with its graph, we obtain an idempotent
$$\tpi_j:=(id-\sigma_j)/2\in\Cor^0_S(D_j,D_j)=
\CH^0(D_j\times_SD_j)_{\Q}$$
in the group of relative correspondences (see 1.6 for the definition).
We define a Chow motive, called the Prym motive, by
$$\Prym(D_j/C'_j):=(D_j,\tpi_j).$$
This is a relative Chow motive, and can be viewed as an absolute
Chow motive, see (1.6.1).
Let $h^i(X)$, $h^i(S)$, $h^i(C'_j)$ denote the $i$-th component of
the Chow-K\"unneth decomposition [17], [18] (their existence was
proved there for $h^i(S)$, $h^i(C'_j)$ and in [1] for $h^i(X)$).
Set $h(X)=\mopl_i\,h^i(X)\,(=(X,\Delta)$ where $\Delta$ is the
diagonal), and similarly for $h(S)$, $h(C'_j)$.

For $j>r$, define as absolute Chow motives
$$\Prym^i(D_j/C'_j):=\Prym(D_j/C'_j)\,\,\,\hbox{if}\,\,\,i=1,\,\,\,
\hbox{and $0$ otherwise}.$$
Then $\Prym^1(D_j/C'_j)$ is identified with the Prym variety of
$D_j/C'_j$ by Weil's theory of correspondences between curves.

If $j\le r$, choosing $\xi_j\in\CH^1(C'_j)_{\Q}$ such that the degree
of its restriction to each irreducible component of
$C'_j\otimes_k\ok$ is 1, we can construct a decomposition as absolute
Chow motives (see (1.11) below)
$$\Prym(D_j/C'_j)=\mopl_{i=0}^2\,\Prym^i(D_j/C'_j),$$
such that we have in case $k=\ok$
$$\Prym^i(D_j/C'_j)\cong h^i(C'_j).$$
However, it does not seem that the last isomorphisms hold in case
$k\ne\ok$, see (1.12).

Let $\ell$ be a prime different from the characteristic of $k$, and
let $\CH^p_{\alg}(X)_{\Q}$ be the subgroup of $\CH^p(X)_{\Q}$
consisting of cycles algebraically equivalent to zero.
The following gives a generalization of [3], [6] (and [1], where the authors proved the existence of a Chow--K\"unneth decomposition), and has been conjectured by the first author [19].

\medskip\noindent
{\bf Theorem~1.} {\it There is a self-dual Chow-K\"unneth
decomposition for $X$ together with isomorphisms of Chow motives
$$h^i(X)\cong h^i(S)\oplus h^{i-2}(S)(-1)\oplus(\mopl_j\,\Prym^{i-2}
(D_j/C'_j)(-1)),$$
where $(-1)$ denotes the Tate twist of Chow motives.
In particular, if $H^1(S_{\ok},\Q_{\ell})=0$ or equivalently
$\CH^1_{\alg}(S_{\ok})_{\Q}=0$, then}
$$h^3(X)\cong\mopl_j\,\Prym^1(D_j/C'_j)(-1).$$

\medskip
Note that if $k=\ok$ or more generally $r=r'$, then the first
isomorphisms become
$$\eqalign{
h^3(X)&\cong h^3(S)\oplus h^{1}(S)(-1)\oplus(\mopl_{j\le r}\,
h^1(C'_j)(-1))\oplus(\mopl_{j>r}\,\Prym(D_j/C'_j)(-1)),\cr
h^i(X)&\cong h^i(S)\oplus h^{i-2}(S)(-1)\oplus(\mopl_{j\le r}\,
h^{i-2}(C'_j)(-1))\hskip2.5cm\hbox{if}\,\,\, i\ne 3.\cr}$$

Theorem~1 gives a refinement (up to an isogeny) of a theorem of
Beauville [3] in the case of conic bundles over $\P^2_{\C}$ with
smooth $C$, where he gave an isomorphism between the intermediate
Jacobian of $X$ and the Prym variety $\cP_X$ of $\tC/C$ as principally
polarized abelian varieties over $\C$.
Note that Theorem~1 in the case $k=\C$ implies an isomorphism of
$\Q$-Hodge structures
$$H^3(X)=H^3(S)\oplus H^1(S)(-1)\oplus(\mopl_j\hbox{Coker}
(H^1(C'_j)\to H^1(D_j))(-1)).$$
To show Theorem~1, we consider the relative Chow-K\"unneth
decomposition for $f$ (see [9], [14], [15], [22]) in the 'weak' and
'strong' sense (see 1.7 for notation), and prove the following
(which has been studied in [19]).

\medskip\noindent
{\bf Theorem~2.} {\it There is a canonical self-dual relative
Chow-K\"unneth decomposition for $f$ in the weak sense, and the
projectors $\pi_{f,-1}$, $\pi_{f,0}$ and $\pi_{f,1}$ define relative
Chow motives isomorphic to
$(S,\Delta_S)$, $\mopl_j\,\Prym(D_j/C'_j)(-1)$ and
$(S,\Delta_S)(-1)$ respectively, where $\Delta_S$ is the diagonal of
$S\times S$.
Moreover, there is a canonical self-dual relative Chow-K\"unneth
decomposition for $f$ in the strong sense, and the relative projector
$\pi_{f,0,j}$ corresponding to the direct factor supported on $C_j$
defines a relative Chow motive isomorphic to $\Prym(D_j/C'_j)(-1)$.
}

\medskip
The proof of Theorem~2 follows from a calculation of the composition
of certain relative correspondences by decomposing these into the
compositions of more elementary correspondences.
Here we have to show the vanishing of certain `phantom' motives.
The construction of the middle projector is due to the first author
[19].
We have the uniqueness of the self-dual decompositions in case $r=0$,
see Remark~(2.6).

\medskip
From Theorem~2 we can deduce the following generalization of [3],
Th.~3.6 (where $k=\ok$ and $S=\P^2$) and [6], Th.~2.6 (where $k=\ok$,
$\char k=0$ and $C$ is irreducible).

\medskip\noindent
{\bf Corollary~1.} {\it There is a canonical isomorphism
$$\CH_{\alg}^2(X)_{\Q}=\CH_{\alg}^2(S)_{\Q}\oplus\CH_{\alg}^1(S)_{\Q}
\oplus\cP_X(k)_{\Q}.$$
In particular, if $H^1(S_{\ok},\Q_{\ell})=0$ or equivalently
$\CH^1_{\alg}(S_{\ok})_{\Q}=0$, then
$$\CH_{\alg}^2(X)_{\Q}=\CH_{\alg}^2(S)_{\Q}\oplus\cP_X(k)_{\Q}.$$
If furthermore $\CH^2(S)_{\Q}=\Q$, then}
$$\CH_{\alg}^2(X)_{\Q}=\cP_X(k)_{\Q}.$$

\medskip
In case $k=\ok$ and $\char k=0$, the condition $\CH^2(S)_{\Q}=\Q$
implies $H^i(S,\O_S)=0$ for $i=1,2$, see [16].
Its converse was conjectured by S.~Bloch [7], and has been proved
at least if $S$ is not of general type, see [8] and also [2], etc.

\medskip
In Section 1 we review some basic facts related to conic bundles and
Chow-K\"unneth decompositions.
In Section 2 we prove the main theorems.

\medskip
This paper grew out of several discussions between the authors.
We would like to thank J.~Murre for useful discussions, and for giving
us the opportunity of the discussions.

\bigskip\bigskip
\centerline{\bf 1.\ Preliminaries}

\bigskip\noindent
{\bf 1.1.~Conic bundles.}
Let $f:X\to S$ be a conic bundle with $\dim X=3$ and $\dim S=2$.
Let $C$ be the discriminant. It is a divisor with normal crossings,
see [3]. Locally $X$ is a subvariety of $U\times\P^2$ defined
by a relative quadratic form where $U$ is an open subvariety of $S$.
Note that $X_s:=f^{-1}(s)$ is a union of two lines (resp.\ a line) in
$\P^2$ if $s$ is a smooth (resp.\ singular) point of $C$.
Let $X_C=f^{-1}(C)$, and let $\tXC$ be its normalization.
Let $C^{\prime}$ be the normalization of $C$. Then $\tXC$ is smooth,
and is a $\P^1$-bundle over a double covering $D$ of $C^{\prime}$
(its fibers are lines in $\P^2$ locally).

Let $C_j$ be the irreducible components of $C$. Let $C_j^{\prime}$
be the normalization of $C_j$, and $D_j$ be the union of the
irreducible components of $D$ whose image in $C$ is $C_j$.
Put $C_j^o=C_j\setminus\Sing C$. In the sequel we shall identify
$C_j^o$ with the corresponding subset of the normalization
$C_j^{\prime}$. Let
$$\rho_j: D_j\to C_j^{\prime}$$
be the double covering, and put $D_j^o = \rho_j^{-1}(C_j^o)$.
Consider the condition:
$$\hbox{The double covering $D_j^o\to C_j^o$ is trivial.}
\leqno(1.1.1)$$
Renumbering the $C_j$ if necessary, there are integers $r\ge r'\ge 0$
such that

\smallskip
Condition (1.1.1) holds if and only if $1\le j\le r'$,

\smallskip
Condition (1.1.1) holds after the base change $k\to\ok$
if and only if $1\le j\le r$.

\medskip\noindent
{\bf 1.2.~Remark.}
In case $C_j^o\otimes_k\ok$ is not connected, condition
(1.1.1) depends only on the restriction over any connected component
of $C_j^o\otimes_k\ok$.
Indeed, let $K_j$ and $K'_j$ denote respectively the function
fields of $C_j^o$ and $D_j^o$, and set $k_j=K_j\cap\ok$,
$k'_j=K'_j\cap\ok$.
Let $n_j$ and $n'_j$ denote respectively the numbers of connected
components of $C_j^o\otimes_k\ok$ and $D_j^o\otimes_k\ok$. Then
$$n_j=\deg k_j/k,\quad n'_j=\deg k'_j/k.\leqno(1.2.1)$$
Hence $\deg k'_j/k_j$ is either 2 or 1, depending on whether the
covering is trivial or not.

To show (1.2.1), let $k''\supset k$ be a sufficiently large finite
Galois extension in $\ok\subset\overline{K'_j}$ containing $k_j$,
$k'_j$.
Then $K_jk''$ is a Galois extension over $K_j$ such that the
restriction induces an isomorphism of Galois groups
$$\Gal(K_jk''/K_j)\simto\Gal(k''/k_j).$$
So $\deg k''/k_j=\deg K_jk''/K_j$, and hence
$$K_j\otimes_{k_j}k''=K_jk'',$$
i.e.\ $C_j^o$ is absolutely irreducible over $k_j$,
where $k_j\subset\Gamma(C_j^o,\O)$ since $C_j^o$ is normal.
(A similar assertion holds for $K'_j$.)
Thus (1.2.1) is proved.

\medskip\noindent
{\bf 1.3.~Example.}
Let $E_j$ be line bundles on $S$, and $a_j$ be sections of
$E_j\otimes E_j$ for $j=0,1,2$. Assume the zeros of $a_j$ are
smooth divisors $C_j$ and their union $C$ is a divisor with normal
crossings on $S$.
Then these define a conic bundle $f:X\to S$ such that $X$ is locally
defined by $$\msum_{0\le j\le 2}\,a_jx_j^2=0\quad\hbox{in}\,\,\,
U\times\P^2,$$
trivializing $E_j$ locally over an open subvariety $U$ of $S$.
Condition (1.1.1) does not hold if $C_j\cap\Sing C\ne\emptyset$,
see (1.4.4) below.

\medskip\noindent
{\bf 1.4.~Decomposition theorem.}
Assume that $k$ is an algebraically closed field.

The decomposition theorem [5] states that if $f:X\to S$ is a proper morphism of irreducible varieties over $k$ such that $X$ is smooth, there is a non canonical
isomorphism ($n = \dim X$)
$$
Rf_*\Q_{\ell,X}[n]\simeq\mopl_i {}^pR^if_*(\Q_{\ell,X}[n])[-j]
$$
in $D^b_c(S,\Q_{\ell})$ together with canonical isomorphisms
$$
^pR^if_*(\Q_{\ell,X}[n])\simeq\mopl_Z \IC_Z(E^i_{Z^o}),
$$
where $Z$ runs over the integral closed subvarieties of $S$, $E^i_{Z^o}$ is a smooth $\Q_{\ell}$-sheaf over a dense open subset $Z^o\subset Z$ and $\{Z\mid E^i_{Z^o}\ne 0\}$ is a finite set. See [5] for the definition of $D_c^b(S,\Q_{\ell})$ and
${}^pR^if_*:={}^p\cH^i\R f_*$.

With the notation and the assumptions of (1.1),
let $\iota_j:C_j^o:=C_j\setminus\Sing C\to C_j$ denote the inclusion.
In our case the decomposition theorem implies the existence of a noncanonical isomorphism
$$\R f_*\Q_{\ell,X}[3]\simeq\mopl_{-1\le i\le 1}{}^pR^if_*(\Q_{\ell,X}
[3])[-i]\quad\hbox{in}\,\,\, D_c^b(S,\Q_{\ell}),\leqno(1.4.1)$$
together with canonical isomorphisms
$$\eqalign{{}^pR^{-1}f_*(\Q_{\ell,X}[3])&=\Q_{\ell,S}[2],\quad
{}^pR^{1}f_*(\Q_{\ell,X}[3])=\Q_{\ell,S}(-1)[2],\cr
{}^pR^{0}f_*(\Q_{\ell,X}[3])&=\mopl_j(\iota_j)_*L_j[1].\cr}
\leqno(1.4.2)$$
Here $L_j$ is the restriction to $C_j^o$ of $((\rho_j)_*\Q_{\ell,D_j}/
\Q_{\ell,C_j^{\prime}})(-1)$ with $\rho_j:D_j\to C_j^{\prime}$ the
natural morphism.
It is a smooth $\Q_{\ell}$-sheaf of rank 1.

Note that
$$\hbox{Condition (1.1.1) is equivalent to
$\Gamma(C_j^o,L_j)\ne 0$.}\leqno(1.4.3)$$
Since the fiber of $f$ at $s\in C_j\setminus C_j^o$ is
a line, the stalk of $(\iota_j)_*L_j$ at $s\in C_j\setminus C_j^o$
vanishes and hence $(\iota_j)_*L_j=(\iota_j)_!L_j$, i.e. the local
monodromy of $L_j$ around $s$ is nontrivial. So we get
$$\hbox{Condition (1.1.1) does not hold if
$C_j\cap\Sing C\ne\emptyset$.}\leqno(1.4.4)$$
Note that the last condition is equivalent to that any connected
component of $C$ has a singular point.

\medskip\noindent
{\bf 1.5.~Chow motives.}
Let $X,Y$ be smooth projective varieties over a perfect field $k$.
Assume $X$ is equidimensional. Then the group of correspondences
is defined by
$$\Cor_k^i(X,Y)=\CH^{\dim X+i}(X\times_kY)_{\Q}.$$
In general, we take the direct sum over the connected components of
$X$.
A {\it Chow motive} is defined by $(X,\pi,i)$ where
$\pi\in\Cor_k^0(X,X)$ is an idempotent (i.e. $\pi^2=\pi$) and $i\in\Z$.
Note that $i$ is related to morphisms of Chow motives
which are defined by
$$\Hom((X,\pi,i),(Y,\pi',j))=\pi'\scirc\Cor_k^{j-i}(X,Y)\scirc\pi.$$
Sometimes we denote $(X,\pi,0)$ by $(X,\pi)$.
The {\it Tate twist} of Chow motives is defined by
$$(X,\pi,i)(m)=(X,\pi,i+m).$$

Similarly we can define relative Chow motives (see [9], [12])
using relative correspondences defined as below.

\medskip\noindent
{\bf 1.6.~Relative correspondences.}
Let $X,Y$ be smooth varieties over a perfect field $k$
with projective morphisms $f:X\to S$, $g:Y\to S$ over $k$.
The group of relative correspondences is defined by
$$\Cor_S^i(X,Y)=\CH_{\dim Y-i}(X\times_SY)_{\Q},$$
if $Y$ is equidimensional. In general we take the direct sum over
the connected components of $Y$.
The composition of relative correspondences is defined by using
the pull-back associated to the cartesian diagram
$$\matrix{X\times_SY\times_SZ&\to&(X\times_SY)\times_k(Y\times_SZ)\cr
\downarrow&&\downarrow\cr Y&\to&Y\times_kY,\cr}$$
together with the pushforward by $X\times_SY\times_SZ\to X\times_SZ$,
see [9], [13].
There is a natural morphism
$$\Cor_S^i(X,Y)\to\Cor_k^i(X,Y),\leqno(1.6.1)$$
which is compatible with composition.
This induces a forgetful functor from the category of relative Chow
motives over $S$ to the category of Chow motives over $k$, see [9].

If $k=\ok$ we have the action of correspondences
$$\eqalign{\Cor_S^i(X,Y)
&\to\Hom(\R f_*\Q_{\ell,X},\R g_*\Q_{\ell,Y}(i)[2i])\cr
&\to\mopl_j\Hom({}^pR^jf_*\Q_{\ell,X},{}^pR^{j+2i}g_*\Q_{\ell,Y}(i)).
\cr}\leqno(1.6.2)$$
This is compatible with the composition of correspondences,
see loc.~cit.

\medskip\noindent
{\bf 1.7.~Relative Chow-K\"unneth decomposition.}
Let $f:X\to S$ be a proper morphism  of irreducible varieties over $k$. We say that $f$ admits a {\it relative Chow-K\"unneth decomposition in the weak sense} if there exist mutually orthogonal idempotents
$\pi_{f,i}\in\Cor_S^0(X,X)$
such that $\sum_i\pi_{f,i}=\Delta_X$ (where $\Delta_X$ denotes the
diagonal) and such that the action of $\pi_{f,i}$ on
$^pR^jf_*(\Q_{\ell,X}[n])$ is the identify for $i=j$, and vanishes
otherwise; see [22].
In case $k$ is not algebraically closed, we say that mutually
orthogonal idempotents define a relative Chow-K\"unneth decomposition
if their base changes by $k\to\ok$ do.

We say that $f:X\to S$ admits a {\it relative Chow-K\"unneth decomposition in the strong sense} if there exist mutually orthogonal idempotents $\pi_{i,Z}\in\Cor^0_S(X,X)$ such that $\sum_{i,Z}\pi_{i,Z}=\Delta_X$ and such that the action of $\pi_{i,Z}$ on $\IC_W(E^j_{W^0})$ is the identity if $(i,Z) = (j,W)$ and zero otherwise.

Note that $f$ admits a relative Chow-K\"unneth decomposition in the strong sense if and only if
$f:X\to S$ satisfies the {\it motivic decomposition conjecture} [9], [14].

In the case of a conic bundle over a surface (notation and assumptions as in (1.4)), assume there are
mutually orthogonal idempotents
$$\pi_{f,i}\in\Cor_S^0(X,X)=\CH^1(X\times_SX)_{\Q}\quad\hbox{for}
\,\,\,i=-1,0,1,$$ defining
a relative Chow-K\"unneth decomposition for $f$, and let $\pi_{f,0,j}$ be
mutually orthogonal relative projectors such that $\pi_{f,0}=
\sum_j\pi_{f,0,j}$. 

(Note that $\pi_{f,i}\scirc\pi_{f,0,j}=\pi_{f,i}
\scirc\pi_{f,0}\scirc\pi_{f,0,j}=0$ for $i=\pm 1$.)
The projectors $\pi_{f,0,j}$ define a relative Chow-K\"unneth decomposition
for $f$ in the strong sense if the action of $\pi_{f,0,j}$ on the
direct factor supported on $C_{j'}$ is the identify for $j=j'$,
and vanishes otherwise.
In case $k$ is not algebraically closed, the above condition should
be satisfied for the base change by $k\to\ok$, where the direct factor
supported on $C_{j'}$ should be replaced by the direct factor supported
on the base change of $C_{j'}$.

We say that a decomposition is {\it self-dual} if the projectors
satisfy the self-duality
$$\pi_{f,i}={}^t\pi_{f,-i}\,\,\,\hbox{and}\,\,\,
\pi_{f,0,j}={}^t\pi_{f,0,j}\,\,\hbox{(in the strong case).}$$

\medskip\noindent
{\bf 1.8.~Heuristic argument.}
With the notation and the assumptions of (1.4), assume that the
decomposition (1.4.1) holds in the derived category of (conjectural)
motivic sheaves $D^b\M(S)$ (see [4]) where the following isomorphism
should hold:
$$\End_{D^b\M(S)}(\R f_*\Q^{\M}_X[3])=\Cor_S^0(X,X)\,
(:=\CH^1(X\times_SX)_{\Q}).\leqno(1.8.1)$$
Here $\Q^{\M}_X\in D^b\M(X)$ is the constant sheaf.
(In case $k=\C$ we may assume $\M(X)=\MHM(X)$, see Remark~(1.9)
below.)
Then (1.4.1) and (1.8.1) should induce a relative Chow-K\"unneth
decomposition in the weak sense by taking the projection to each
direct factor.
If we have another relative Chow-K\"unneth decomposition in the weak
sense, then the corresponding projectors $\pi_{f,i}$ are identified
with endomorphisms
$$\pi_{f,i}:\R f_*\Q^{\M}_X[3]\to\R f_*\Q^{\M}_X[3],$$
and (1.4.1) gives a decomposition
$\pi_{f,i}=\mopl_{a,b}(\pi_{f,i})_{a,b}$ such that $(\pi_{f,i})_{a,b}$
is identified with
$$(\pi_{f,i})_{a,b}\in\Ext^{a-b}({}^pR^bf_*(\Q^{\M}_X[3]),
{}^pR^af_*(\Q^{\M}_X[3]))\,\,\,\,(i,a,b\in\{-1,0,1\}).$$
In particular, $(\pi_{f,i})_{a,b}=0$ for $a>b$. We have also
$$\hbox{$(\pi_{f,i})_{i,i}=id$, and $(\pi_{f,i})_{a,a}=0$ for $i\ne a
\,\,(i,a\in\{-1,0,1\})$.}$$

Assume now $r=0$, i.e.\ (1.1.1) does not not hold for any $j$. Then
$$(\pi_{f,i})_{a,b}=0\quad\hbox{if}\,\,\, a-b=1.\leqno(1.8.2)$$
Indeed, for $(a,b)=(0,1)$ we have by (1.4.3)
$$\Ext^1(\Q_{\ell,S}(-1)[2],\mopl_j(\iota_j)_*L_j[1])=
\mopl_jH^0(C_j^o,L_j)(1)=0.$$
For $(a,b)=(-1,0)$, the assertion follows from duality
since $L_j(1)$ is self-dual.

By (1.8.2) we have for $i=-1,0,1$
$$\pi_{f,i}=(\pi_{f,i})_{i,i}+(\pi_{f,i})_{-1,1}.$$
It is then easy to see that the condition
$\pi_{f,0}\scirc\pi_{f,0}=\pi_{f,0}$ implies
$$\hbox{$(\pi_{f,0})_{-1,1}=0$, i.e. $\pi_{f,0}=(\pi_{f,0})_{0,0}$.}$$
In particular, $\pi_{f,0}$ is {\it unique} in this case. Note that
$(\pi_{f,1})_{-1,1}+(\pi_{f,-1})_{-1,1}=0$ by $\pi_{f,-1}\scirc
\pi_{f,1}=0$, and $(\pi_{f,i})_{-1,1}$ for $|i|=1$ gives the ambiguity
of the decomposition. Indeed, for any $\eta\in\Ext^2(\Q_{\ell,S}(-1)
[2],\Q_{\ell,S}[2])$, we can replace $\pi_{f,1}$, $\pi_{f,-1}$ with
$\pi_{f,1}+\eta$ and $\pi_{f,-1}-\eta$ respectively.
(If we assume the self-duality of the decomposition,
this imposes some condition on the ambiguity.)

If $r>0$, then (1.8.2) does not hold, and the situation is rather
complicated.
It is not clear whether the uniqueness of the decomposition holds
even the self-duality is assumed.

\medskip\noindent
{\bf 1.9.~Remark.}
In case the base field is $\C$, the above argument can be
justified. Indeed, let $d_X=\dim X$ and $Y=X\times_SX$ with
the projections $pr_i:Y\to X$. Let $\D_Y$ denote the dualizing
complex. Then, using the adjunction and the base change in [20],
we have the isomorphisms (see also [9])
$$\eqalign{\End_{D^b\MHM(S)}(\R f_*\Q_X)
&=\Hom_{D^b\MHM(X)}(\Q_X,f^!\R f_*\Q_X)\cr
&=\Hom_{D^b\MHM(X)}(\Q_X,\R(pr_1)_*pr_2^!\Q_X)\cr
&=\Hom_{D^b\MHM(Y)}(pr_1^*\Q_X,pr_2^!\Q_X)\cr
&=\Ext^{-2d_X}_{D^b\MHS}(\Q,\R\Gamma(Y,\D_Y(-d_X)).\cr}
$$
Here $\MHS$ and $\MHM(X)$ denote respectively the categories of
polarizable mixed Hodge structures [10] and mixed Hodge modules on
$X$ [20]. We have moreover the following

\medskip\noindent
{\bf 1.10.~Proposition.} {\it
Let $Y$ be a complex algebraic variety such that $\dim\Sing Y\le d_Y-2$
where $d_Y=\dim Y$. Then we have an isomorphism}
$$
\CH^1(Y)_{\Q}=\Ext^{2-2d_Y}_{D^b\MHS}(\Q,\R\Gamma(Y,\D_Y(1-d_Y)).
$$

\medskip\noindent
{\it Proof.}
Let $Z=\Sing Y$ and $U=Y\setminus Z$ with the inclusions $i:Z\to Y$
and $j:U\to Y$. Since $\dim Z\le\dim Y-2$, we have
$$\CH^1(Y)=\CH^1(U).$$
On the other hand, there is a distinguished triangle in $D^b\MHM(Y)$
$$i_*\D_Z\to\D_Y\to\R j_*\D_U\to,$$ inducing a long exact sequence of
extension groups $\Ext^i_{D^b\MHS}(\Q,\R\Gamma(Y,*))$, and
$$\Ext^{-i}_{D^b\MHS}(\Q,\R\Gamma(Z,\D_Z(1-d_Y))=0\quad
\hbox{for}\,\,i>2\dim Z,$$ since
$$H^{\BM}_i(Z)=H^{-i}(Z,\D_Z)=0\quad\hbox{for}\,\,i>2\dim Z.$$
So the assertion is reduced to the smooth case, and follows from
[21], Prop.~3.4.
This finishes the proof of Proposition~(1.10).

\medskip\noindent
{\bf 1.11.~Decomposition of Prym motives.}
Let $\rho:X\to Y$ be a surjective finite morphism of algebraic
varieties over a perfect field $k$.
Assume $X$ is smooth over $k$, $Y$ is irreducible, $\rho$ is
generically of degree 2, and $\char k\ne 2$.
Let $U$ be a non-empty open subvariety of $Y$ over which
$\rho$ is finite \'etale of degree 2. Let
$$\sigma\in\Cor_Y^0(X,X)=\CH_{\dim X}(X\times_YX)_{\Q}
=\CH_{\dim X}(X_U\times_UX_U)_{\Q},$$
such that its restriction over $U$ is the involution associated to
the finite \'etale covering of degree 2, where $X_U=\rho^{-1}(U)$.
The relative Prym motive is defined by
$$\Prym(X/Y)=(X,\pi)\quad\hbox{with}\quad\pi=(id-\sigma)/2.
\leqno(1.11.1)$$

Assume now that $X=X_1\coprod X_2$.
There is a canonical isomorphism
$$X_1\cong X_2\quad\hbox{over}\,\,\,Y,\leqno(1.11.2)$$
since $X_1,X_2$ are the normalization of $Y$ (because they
are smooth over $k$ and finite over $Y$).
For $a=1,2$, we have an isomorphism
$$\Prym(X/Y)\cong(X_a,\Delta).\leqno(1.11.3)$$
Indeed, using the restriction over $U$, we get
$$\Cor_Y^0(X_a,X)=\Cor_Y^0(X,X_a)=\Q\oplus\Q,$$
and (1.11.3) is reduced to
$$\eqalign{\pmatrix{1\cr -1\cr}\pmatrix{1/2&-1/2}&=
\pmatrix{1/2&-1/2\cr-1/2&1/2\cr}\cr\pmatrix{1/2&-1/2}
\pmatrix{1\cr -1\cr}&=1.\cr}$$
Here $\pi$ is represented by $\pmatrix{1/2&-1/2\cr-1/2&1/2\cr}$
using (1.11.2), and the isomorphism (1.11.3) is induced by
$\pmatrix{c\cr -c\cr}$ and $\pmatrix{1/2c&-1/2c}$ which are
identified respectively with elements of $\Cor_Y^0(X_a,X)$ and
$\Cor_Y^0(X,X_a)$, where $c$ is any nonzero rational number.

If moreover $X,Y$ are projective, then choosing a 0-cycle
$\xi_a\in\CH_0(X_a)_{\Q}$ with degree 1 in a compatible way with
(1.11.2), the 0-th K\"unneth projector $\pi_{X_a,0}$ of
$(X_a,\Delta)$ is defined by
$$\xi_a\times[X_a].$$
Composing it with $\pmatrix{c\cr -c\cr}$ and $\pmatrix{1/2c&-1/2c}$,
we get the projector defining $\Prym^0(X/Y)$, and it is explicitly
expressed by
$$\sum_{a=1}^2\xi_a\times[X_a]/2-\sum_{a=1}^2\xi_a\times[X_{3-a}]/2
\in\Cor_k^0(X,X).$$

Assume now that $X$ is projective and irreducible, but
$$k'':=k(X)\cap\ok\ne k(Y)\cap\ok=:k'.$$
Here we choose an embedding $\ok\to\overline{k(X)}$. It induces
$k(X)\otimes_k\ok\to\overline{k(X)}$ and defines a geometric generic
point of an irreducible component $X_{\ok,0}$ of
$X_{\ok}:=X\otimes_k\ok$.
Let $G'$ and $G''$ denote the subgroups of the Galois group $G$ of
$\ok/k$ corresponding to $k'$ and $k''$ respectively.
Set $d'=|G/G'|$, $d''=|G/G''|$ so that $d''=2d'$.
Let $\xi\in\CH_0(X)_{\Q}$ with degree $d''$ so that its restriction
$\xi_0$ to $X_{\ok,0}$ has degree 1.
Take $g_i\in G\,(i\in[1,d'])$, $h\in G'$ such that
$$\hbox{$G=\coprod_{i=1}^{d'}g_iG',\quad G'=G''\coprod hG''$,\quad
hence\quad$G=\coprod_{i=1}^{d'}(g_iG''\coprod g_ihG'')$}.$$
Here $h^2G''=G''$ since $|G'/G''|=2$. (So $g_i$ can be replaced by
$g_ih$.) The projector defining $\Prym^0(X/Y)$ is then given by
$$\sum_{i=1}^{d'}(g_i\xi_0\times g_i[X_{\ok,0}]-
g_i\xi_0\times g_ih[X_{\ok,0}]-g_ih\xi_0\times g_i[X_{\ok,0}]+
g_ih\xi_0\times g_ih[X_{\ok,0}])/2.$$
This is invariant by the action of $G$ since $g\xi_0=\xi_0$,
$g[X_{\ok,0}]=[X_{\ok,0}]$ for $g\in G''$.
The argument is similar for $\Prym^{2n}(X/Y)$ with $n=\dim X$
(exchanging the first and second factors of the product).

\medskip\noindent
{\bf 1.12.~Remark.}
Let $C$ be a smooth projective curve over $k$, and $D=C\otimes_kk'$
where $k'/k$ is a field extension of degree 2 (and $\char k\ne 2$).
Assume $k'\not\subset k(C)$ so that $D$ is irreducible.
Then
$$\Cor_C^0(D,D)=\CH^0(D\times_CD)_{\Q}=\Q\oplus\Q,$$
since 
$$D\times_CD=C\otimes_k(k'\otimes_kk')=C\otimes_k(k'\oplus k')=
\hbox{$D\coprod D$},$$
where the two $D$ correspond to the diagonal and the antidiagonal.
So we can define the Prym motive $\Prym(D/C)$ as in (1.11.1).
However, it does not seem that
$$\Prym(D/C)\cong(C,\Delta),\leqno(1.12.1)$$
without taking the base change $k\to\ok$.
It does not hold at least over $C$, since
$$\Cor_C^0(D,C)=\CH^0(D)_{\Q}=\Q.$$
The problem is whether they are isomorphic over $k$, and we have to
consider
$$\Cor_k^0(D,C)=\CH^1((C\times_kC)\otimes_kk')_{\Q}=
\CH^1(D\times_{k'}D)_{\Q}.$$
In the case $C$ is an elliptic curve without complex multiplication,
it does not seem that the above group contains an element inducing
the desired isomorphism.

Related to this, we have the following problem:
$$\hbox{Is there an isomorphism}\,\,\,\,
\Prym(k'/k)\cong(\Spec k,\Delta)\,?\leqno(1.12.2)$$
Here the left-hand side is defined as in (1.11.1).
Note that (1.12.2) holds after taking the base change $k\to\ok$.
However, it does not hold without the base change since
$$\Cor_k^0(\Spec k,\Spec k')=\CH^0(\Spec k')_{\Q}=\Q.$$
Note also that (1.12.1) should imply (1.12.2) in case
$\ok\cap k(C)=k$ since we should have
$$\Prym^0(D/C)\cong\Prym(k'/k),\quad h^0(C)\cong(\Spec k,\Delta).
\leqno(1.12.3)$$ 

\bigskip\bigskip
\centerline{\bf 2.\ Proof of main theorems}

\bigskip\noindent
{\bf 2.1.~Lemma.} {\it With the notation of $(1.6)$, assume $f,g$
are flat. Set $n=\dim X-\dim S$.
Let $\xi\in\Cor_S^i(X,S)=\CH^{n+i}(X)_{\Q}$ and
$\xi'\in\Cor_S^j(S,Y)=\CH^{j}(Y)_{\Q}$.
Let $pr_1:X\times_SY\to X$ and $pr_2:X\times_SY\to Y$ denote the
projections.
Then the composition $\xi'\scirc\xi\in\Cor_S^{i+j}(X,Y)=
\CH^{n+i+j}(X\times_S Y)_{\Q}$ is given by $pr_1^*\xi$ if $\xi'=[Y]$,
and $pr_2^*\xi'$ if $\xi=[X]$.
}

\medskip\noindent
{\it Proof.} The flatness of $f,g$ implies that $X\times_SY\to
X\times_kY$ is a regular embedding and the $pr_i$ are flat.
Moreover, we have locally a regular sequence defining
$X\times_SY$ in $X\times_kY$ and it is a regular sequence for the
pull-back by $X\times_kY\to Y$ of any $\O_Y$-module.
The last assertion follows from the flatness of $pr_2$ together with
the theory of regular sequences (see e.g. [23], p.~71)
since the Koszul complex calculates the pull-back by the embedding
$X\times_SY\to X\times_kY$. So the assertion follows.

\medskip\noindent
{\bf 2.2.~Lemma.} {\it With the notation of $(1.6)$, let
$\xi\in\Cor^i_S(S,X)=\CH^i(X)_{\Q}$ and
$\xi'\in\Cor^j_S(X,S)=\CH^{j+n}(X)_{\Q}$ where $n=\dim X-\dim S$.
Then $\xi'\scirc\xi\in\Cor_S^{i+j}(S,S)=\CH^{i+j}(S)_{\Q}$ is given
by $f_*(\xi\cdot\xi')\in\CH^{i+j}(S)_{\Q}$, where $\xi\cdot\xi'$
is the intersection of cycles on $X$.
}

\medskip\noindent
{\it Proof.} This immediately follows from the definition of the
composition in (1.6).

\medskip\noindent
{\bf 2.3.~Lemma.} {\it With the notation of $(1.1)$, let
$\xi,\xi'\in\Cor^0_S(X,X)=\CH^1(X\times_SX)_{\Q}$ which are
represented by cycles supported in the inverse images of curves $C$
and $C'$ respectively on $S$.
Assume $\dim C\cap C'\le\dim S-2$ or one of the cycles belongs to
$pr^*\CH^1(S)_{\Q}$ where $pr:X\times_SX\to S$ is the projection.
Then their composition vanishes.
}

\medskip\noindent
{\it Proof.} If the second assumption is satisfied,
we may assume that $\dim C\cap C'\le\dim S-2$ by the
moving lemma on $S$, since one of the cycles comes from $S$. Then the
composition in $\CH^1(X\times_SX)_{\Q}$ is represented by a cycle
supported in the inverse image of $C\cap C'$ which has codimension 2.
So it vanishes.
This finishes the proof of Lemma~(2.3).

\medskip\noindent
{\bf 2.4.~Proof of Theorem~2.}
We first assume that $k$ is algebraically closed.
Take any $\xi\in\CH^1(X)_{\Q}$ such that $f_*\xi=[S]$, i.e. its
restriction to the generic fiber of $f$ is a zero-cycle of degree 1.
The ambiguity of $\xi$ is given by $f^*\eta$ for $\eta\in
\CH^1(S)_{\Q}$.
If $s\notin\Sing C$, there is an open neighborhood $U$ of $s$ such
that the restriction of $\xi$ over $U$ is represented by $[Z]/2$,
where $Z$ is finite \'etale of degree 2 over $U$ since $f$ is a conic
bundle. Set
$$p=pr_1^*\xi\in\Cor_S^0(X,X)=\CH^1(X\times_SX)_{\Q}\quad\hbox{so
that}\quad{}^tp=pr_2^*\xi,$$ where $pr_i$ is the $i$-th projection.
By Lemma~(2.1), we have
$$p=[X]\scirc\xi.$$
where $\xi\in\Cor_S^0(X,S)=\CH^1(X)_{\Q}$ and $[X]\in\Cor_S^0(S,X)=
\CH^0(X)_{\Q}$.
Then $p$ and $^tp$ are idempotents since we have by Lemma~(2.2)
$$\xi\scirc[X]=id\in\Cor_S^0(S,S).\leqno(2.4.1)$$
We have ${}^tp\scirc p=0$ since $^t[X]\scirc[X]=0$ in
$\Cor_S^{-1}(S,S)=0$.
Note that $p\scirc {}^tp = pr^*\eta$ with $\eta=f_*(\xi\cdot\xi)\in
\CH^1(S)_{\Q}$ by Lemmas~(2.1) and (2.2), where $pr:X\times_SX\to S$
is the projection. So we can define
$$\pi_{f,-1}=p\scirc(1-{}^tp/2),\quad\pi_{f,1}=(1-p/2)\scirc{}^tp.$$
Indeed, setting $\pi_{f,-1}=p\scirc(1-a\,{}^tp)$ and $\pi_{f,1}=(1-b
\,p)\scirc {}^tp$ with $a,b\in\Q$, we get $a+b=1$ from the condition
$\pi_{f,-1}\scirc\pi_{f,1}=0$, and $a=b=1/2$ from the self-duality.
Note that $\pi_{f,-1}$ and $\pi_{f,1}$ are still of the form
$pr_1^*\xi$ and $pr_2^*\xi$ respectively, replacing $\xi$ with
$\xi-f^*\eta/2$.
Moreover, they are well-defined. Indeed, if we replace $\xi$ with
$\xi+f^*\zeta$, then $\eta=f_*(\xi\cdot\xi)$ is replaced by
$\eta+2\zeta$, and hence $\pi_{f,-1}$ and $\pi_{f,1}$ are unchanged.

We get thus canonical isomorphisms of relative Chow motives
$$(X,\pi_{f,-1})=(S,\Delta_S),\quad(X,\pi_{f,1})=(S,\Delta_S)(-1),
\leqno(2.4.2)$$
induced by $\xi\in\Cor_S^0(X,S)$ and $^t[X]\in\Cor_S^{-1}(X,S)$ with
inverse $[X]\in\Cor_S^0(S,X)$ and $^t\xi\in\Cor_S^1(S,X)$ respectively.

The action of $\pi_{f,-1}$ on ${}^pR^jf_*(\Q_{\ell,X}[3])$ is the
identity for $j=-1$ and vanishes otherwise (and similarly for
$\pi_{f,1}$), since we have a factorization
$$(\pi_{f,-1})_*:\R f_*(\Q_{\ell,X}[3])\buildrel{\xi_*}\over\lra
(\Q_{\ell,S}[2])[1]\buildrel{[X]_*}\over\lra\R f_*(\Q_{\ell,X}[3]).$$

We now construct the middle projector $\pi_{f,0}$.
Let $\tX_j\subset\tXC$ be the inverse image of $C_j$ and let
$g_j:\tX_j\to X$, $p_j:\tX_j\to D_j$ be natural morphisms. Set
$$\gamma_j := (p_j)_*\scirc{}(g_j)^*\in\Cor_S^{-1}(X,D_j),\quad
\gamma_j':=-{}^t\gamma_j/2\in\Cor_S^1(D_j,X).$$
Let $\sigma_j$ be the involution of $D_j$ associated with
the double covering $D_j/C'_j$. This is identified with a cycle
defined by its graph. The projector $\pi_{f,0,j}$ corresponding
to $C_j^{\prime}$ is defined as in [19] by
$$\pi_{f,0,j}=\gamma_j'\scirc\tpi_j\scirc\gamma_j\quad\hbox{with}\quad
\tpi_j:=(id-\sigma_j)/2.$$
This is represented by a cycle supported in $pr^{-1}(C_j)$,
but does not belong to $pr^*\CH^1(S)_{\Q}$. More precisely,
$\tX_j\times_S\tX_j$ has two irreducible components corresponding to
the compositions of correspondences
$$(p_j)^*\scirc id\scirc(p_j)_*\,\,\,\hbox{and}\,\,\,(p_j)^*\scirc
\sigma_j\scirc(p_j)_*.$$
Taking the composition with $(g_j)^*$ and $(g_j)_*$, we get
the pushforward of these cycles by $g_j$.

By Proposition (2.5) below, $\gamma_j\scirc{}^t\gamma_j\in\Cor_{S}^0
(D_j,D_j)=\Cor_{C_j^o}^0(D^o_j,D^o_j)$ is expressed by the
matrix
$$A:=\pmatrix{-1&1\cr 1&-1\cr}.$$
Here $D^o_j\to C^o_j$ is the restriction of $\rho_j:D_j\to
C_j^{\prime}$ over $C_j^o$; it is \'etale of degree 2.
On the other hand, $\tpi_j:=(id-\sigma_j)/2$ is expressed by the
matrix $$-{1\over 2}A=\pmatrix{1/2&-1/2\cr -1/2&1/2\cr}$$
which is an idempotent since $A^2=-2A$. We get thus
$$\tpi_j=\gamma_j\scirc\gamma'_j.\leqno(2.4.3)$$
This implies that $\pi_{f,0,j}$ is an idempotent, and moreover
$$\pi_{f,0,j}\scirc\gamma'_j\scirc\tpi_j\scirc\gamma_j\scirc
\pi_{f,0,j}=\pi_{f,0,j},\quad\tpi_j\scirc\gamma_j\scirc\pi_{f,0,j}
\scirc\gamma'_j\scirc\tpi_j=\tpi_j.$$
Thus we get an isomorphism of relative Chow motives over $S$
$$(X,\pi_{f,0,j})=\Prym(D_j/C'_j)(-1).$$
Using the compatibility of (1.6.2) with the composition of
correspondences, we get then
$$(\pi_{f,0,j})_*({}^pR^0f_*(\Q_{\ell,X}[3]))=(\iota_{j})_*L_{j}[1]
\quad\hbox{in}\quad{}^pR^0f_*(\Q_{\ell,X}[3]),$$
i.e. the action of the idempotent $\pi_{f,0,j}$ on
$(\iota_{j'})_*L_{j'}[1]\subset{}^pR^0f_*(\Q_{\ell,X}[3])$ is the
identity if $j=j'$, and vanishes otherwise.
The action of $\pi_{f,0,j}$ on $^pR^if_*(\Q_{\ell,X}[3])$ vanishes for
$|i|=1$, since $\pi_{f,0,j}$ is supported in the inverse image of
$C_j$.
Moreover it follows from Lemma~(2.3) that
$$\pi_{f,0,j}\scirc\pi_{f,0,j'}=0\quad\hbox{for}\,\,\,j\ne j'.$$
So we get the middle projector
$$\pi_{f,0}:=\mopl_j\,\pi_{f,0,j}.$$

We have to show the relation
$$\pi_{f,0,j}\scirc\pi_{f,i}=\pi_{f,i}\scirc\pi_{f,0,j}=0
\quad\hbox{for}\,\,\,|i|=1.\leqno(2.4.4)$$
Here it is enough to show $\pi_{f,i}\scirc\pi_{f,0,j}=0$ by duality.
For $i=-1$, it is reduced to
$$\xi\scirc{}^t\gamma_j\scirc(id-\sigma_j)=0\quad\hbox{in}\,\,\,
\Cor_S^1(D_j,S)=\CH^0(D_j)_{\Q}.$$
It is then enough to show the vanishing of its action on
$\Q_{\ell}$-complexes on $S$
$$(\rho_j)_*\Q_{\ell,D_j}\to(\rho_j)_*\Q_{\ell,D_j}\to
\R f_*\Q_{\ell,X}(1)[2]\to\Q_{\ell,S}(1)[2],$$
where the morphisms are induced by $(id-\sigma_j)$, ${}^t\gamma_j$
and $\xi$. So the assertion follows by taking a general transversal
slice $T$ as in the proof of Proposition~(2.5) below.
(Indeed, $\xi$ is represented by $[Z]/2$ over a dense open subvariety
of $S$ where $Z$ is a general hyperplane section of the $\P^2$-bundle
over $S$ containing the conic bundle $X$,
and the assertion is reduced to the vanishing of the intersection
number of $[Z_T]$ and $[X'_{s}]-[X''_{s}]$ in $\overline{X_T}$.
Here $X_T=f^{-1}(T)$, $Z_T=Z\cap X_T$, $\overline{X_T}$ is a
smooth compactification of $X_T$, and $X'_s,X''_s$ are the
irreducible components of $X_s$ where $\{s\}=T\cap D_j$.)
For $i=1$, the assertion is trivial since
$^t[X]\scirc{}^t\gamma_j\scirc(id-\sigma_j)$ belongs to
$\CH^{-1}(D_j)_{\Q}=0$.

Now we have to show
$$\zeta:=1-\msum_{-1\le i\le 1}\,\pi_{f,i}=0\quad\hbox{in}\,\,\,
\Cor_S^0(X,X).\leqno(2.4.5)$$
It is enough to show that $\zeta$ is nilpotent since it
is an idempotent.
As $\CH^0(pr^{-1}(C_j))_{\Q}$ is $4$-dimensional for $j\le r'$,
and is $2$-dimensional otherwise, we have
$$\zeta = pr^*\eta+\msum_{j\le r'}(pr_1^*\xi_j+pr_2^*\xi'_j)
+\msum_j\,c_j\pi_{f,0,j},$$ where $\eta\in\CH^1(S)_{\Q}$,
$\xi_j,\xi'_j\in\CH^0(f^{-1}(D_j))_{\Q}$ and $c_j\in\Q$.
We have $c_j=0$ considering the action of $\zeta$ on
$(\iota_j)_*L_j[1]\subset{}^pR^0f_*(\Q_{\ell,X}[3])$ which vanishes
by the definition of $\zeta$.
(Indeed, the action of $pr^*\eta$, $pr_1^*\xi_j$, $pr_2^*\xi'_j$ on
it vanishes by the same argument as in the case of $\pi_{f,\pm 1}$
using the factorization $pr_1^*\xi=[X]\scirc\xi$, etc.)
By Lemma~(2.3), the assertion (2.4.5) is then reduced to
$$\eqalign{&pr_1^*\xi_j\scirc pr_1^*\xi_j=0,\quad
pr_2^*\xi'_j\scirc pr_2^*\xi'_j=0,\cr
&pr_2^*\xi'_j\scirc pr_1^*\xi_j=0,\quad
pr_1^*\xi_j\scirc pr_2^*\xi'_j\in pr^*\CH^1(S)_{\Q}.\cr}
$$
Here $pr_1^*\xi_j=[X]\scirc\xi_j$ and
$pr_2^*\xi'_j={}^t\xi'_j\scirc{}^t[X]$ in the notation of (2.4.1).
We have the first vanishing since $\xi_j\scirc[X]=f_*\xi_j=0$ using
Lemma~(2.2), and similarly for the second.
The third vanishing follows from the fact that $^t[X]\scirc[X]$
belongs to $\Cor_S^{-1}(S,S)=0$.
For the last assertion, note that $\xi_j\scirc{}^t\xi'_j\in
\Cor_S^1(S,S)=\CH^1(S)_{\Q}$.
So (2.4.5) follows.

Thus Theorem~2 is proved in the case $k=\ok$.
The assertion in the case $k\ne\ok$ is reduced to the case $k=\ok$
since the construction of the relative Chow-K\"unneth projectors is
compatible with the base change although the decomposition of the
middle projector becomes finer after the base change.
So Theorem~2 follows.

\medskip
To complete the proof of Theorem~2 we have to show the following.
(In case $C$ is smooth and irreducible, this also follows from
[9], Example.~5.18.)

\medskip\noindent
{\bf 2.5.~Proposition.} {\it With the above notation, $\gamma_j\scirc
{}^t\gamma_j\in\Cor_S^0(D_j,D_j)=\Cor_{C_j^o}^0(D_j^o,D_j^o)$
is expressed by the matrix $A$.}

\medskip\noindent
{\it Proof.}
Take a sufficiently general closed point $s$ of $C_j^o$.
For $s'\in D_j$ lying over $s$, let $\tX_{s'}$ denote the irreducible
component of $X_s$ corresponding to $s'$ (this is identified with
$p_j^{-1}(s')\subset\tX_j$).
Let $T$ be a sufficiently general transversal slice to $C_j^o$ at $s$,
which is defined by
$$T=h^{-1}(c)\setminus(C_j\setminus\{s\})\,\,\,\hbox{for a
sufficiently general}\,\, c\in k,\leqno(2.5.1)$$
where $T\cap C_j=\{s\}$ and $h$ is a function defined on a
non-empty open subvariety $U$ of $S$ such that $dh\ne 0$ on $U$ and
$dh|_{U\cap T}\ne 0$ on $U\cap T$.
Let $\overline{X_T}$ be a smooth compactification of $X_T:=f^{-1}(T)$
(this exists since it is 2-dimensional). The intersection matrix of
$\tX_{s'},\tX_{s''}$ in $\overline{X_T}$ (where $s',s''$ are the points
of $D_j$ over $s\in C_j^o$) is given by the matrix $A$ since
$[X_s]\cdot[\tX_{s'}]=0$ in $\overline{X_T}$ where we may assume that
$f_T:X_T\to T$ is extended to $\overline{X_T}\to\overline{T}$.

As we have the injection
$$\Cor_S^0(D_j\,D_j)\subset\End((\rho_j)_*\Q_{\ell}),$$ where
$\rho_j:D_j\to C_j$ is the projection, it suffices to calculate the
composition
$$
(\rho_j)_*\Q_{\ell}\buildrel{^t \gamma_j}\over\rightarrow R^2f_*
\Q_{\ell}(1)\buildrel{\gamma_j}\over\rightarrow (\rho_j)_*\Q_{\ell}.
$$
Here the first morphism naturally factors through
$$^t\gamma_j:(\rho_j)_*\Q_{\ell}\to\cH^2_C\R f_*\Q_{\ell}(1),$$
which is the dual of the last morphism, where $\cH^2_C$ is the local
cohomology sheaf.

Restricting these to the transversal slice $T$, we obtain
$$
(\gamma_j)_T\scirc({}^t\gamma_j)_T:\Q_{\ell,s'}\oplus \Q_{\ell,s''}\to
R^2(f_T)_*\Q_{\ell}(1)\to\Q_{\ell,s'}\oplus \Q_{\ell,s''},
\leqno(2.5.2)
$$
where $f_T:X_T\to T$ is the restriction of $f$ over $T$ and similarly
for $(\gamma_j)_T$, $({}^t\gamma_j)_T$.
Here $\Q_{\ell,s'}\oplus\Q_{\ell,s''}$ is identified with a sheaf
supported on $s$.
The first morphism of (2.5.2) naturally factors through
$$({}^t\gamma_j)_T:\Q_{\ell,s'}\oplus \Q_{\ell,s''}\to\H^2_{\{s\}}
(\R(f_T)_*\Q_{\ell}(1)).$$
By the generic base change theorem ([11], 2.9 and 2.10) this is the
dual of the last morphism of (2.5.2) if $c\in k$ in (2.5.1) is
sufficiently general. We have to show that (2.5.2) is expressed by
the intersection matrix $A$.

For $t$, $u\in\{s',s''\}$, the $(t,u)$-component of (2.5.2) is given
by the composition of morphisms of $\ell$-adic cohomology groups
$$
H^0(\{t\})\buildrel{p_j^*}\over\to
H^0(\tX_{t})\buildrel{(\lambda_t)_*}\over\to
H^2_c(X_T)(1)\to
H^2(X_T)(1)\buildrel{(\lambda_u)^*}\over\to
H^2(\tX_{u})(1)\buildrel{p_{j*}}\over\to
H^0(\{u\}),
$$
where $\lambda_t:\tX_t\to X_T$ is the restriction of $g_j$, and
similarly for $\lambda_u:\tX_u\to X_T$.
This is shown by using the commutative diagram
$$\matrix{\H^2_{\{s\}}(K)&\to&(\cH^2K)_s\cr\downarrow&&\uparrow\cr
\H^2_c(T,K)&\to&\H^2(T,K),\cr}$$
where $K=\R (f_T)_*\Q_{\ell}(1)$ so that $\H^2_c(T,K)=H^2_c(X_T)(1)$,
etc.

Moreover the middle morphism $H^2_c(X_T)(1)\to H^2(X_T)(1)$ naturally
factors through $H^2(\overline{X_T})(1)$, and hence we can replace
$X_T$ with $\overline{X_T}$ in the above composition of morphisms.
This implies that (2.5.2) is expressed by the intersection matrix $A$
as is desired. So Proposition~(2.5) follows.

\medskip
As for the uniqueness of the decomposition, it is rather complicated
if $r>0$. However, for $r=0$ we have the following.

\medskip\noindent
{\bf 2.6.~Proposition.} {\it
If $r=0$, the self-dual relative Chow-K\"unneth decomposition is
unique.}

\medskip\noindent
{\it Proof.} Let $\tpi_{f,i}$ be other
mutually orthogonal projectors whose action on the cohomological
direct images is the same as $\pi_{f,i}$. Then $\tpi_{f,i}=\pi_{f,i}$
over a sufficiently small open subvariety of $S$. Hence we have
by the same argument as above (using the condition $r=0$)
$$\tpi_{f,i}=\pi_{f,i}+pr^*\eta_i+\msum_j\,a_{i,j}\pi_{f,0,j}\quad
\hbox{with}\,\,\eta_i\in\CH^1(S)_{\Q},\,\,a_{i,j}\in \Q.$$ We have
$a_{i,j}=0$ by looking at the action on $^pR^0f_*(\Q_{\ell,X}[3])$.
We also get $\eta_0=0$ by $\tpi_{f,0}\scirc\tpi_{f,0}=\tpi_{f,0}$
together with Lemma~(2.3). Moreover, $\eta_{-1}+\eta_1=0$ by
$\tpi_{f,-1}\scirc\tpi_{f,1}=0$ since
$$pr^*\eta_{-1}\scirc\pi_{f,1}=pr^*\eta_{-1},\quad
\pi_{f,-1}\scirc pr^*\eta_1=pr^*\eta_1.$$
(Indeed, for $\xi_1\in\Cor_S^1(S,X)=\CH^1(X)_{\Q}$ and
$\xi_2\in\Cor_S^0(X,S)=\CH^1(X)_{\Q}$, we have $\xi_2\scirc\xi_1=
f_*(\xi_1\cdot\xi_2)\in \CH^1(S)_{\Q}$ by Lemma~(2.2), and this is
$\eta$ in case $\xi_1=\xi$ and $\xi_2=f^*\eta$ since we can take
a good representative of $\xi$ as remarked at the beginning of
this subsection. So the above equalities follow from Lemma~(2.1).)
Then the self-duality implies $\eta_{-1}=\eta_1=0$, and the
uniqueness of the decomposition follows.

\medskip\noindent
{\bf 2.7.~Proof of Theorem~1.}
With the notation of (2.4), we have
$$\pi_{f,-1}=[X]\scirc\xi,\quad\pi_{f,1}={}^t\xi\scirc{}^t[X].$$
Let $\pi_{S,i}$ be the Chow-K\"unneth decomposition for $S$ in [17]
where $\pi_{S,i}=0$ for $i\notin[0,4]$. We may assume the self-duality
$\pi_{S,i}={}^t\pi_{S,4-i}$ as is well-known (by the same argument as
in the construction of $\pi_{f,\pm 1}$ in (2.4)). Define
$$\pi_{X,i}=[X]\scirc\pi_{S,i}\scirc\xi+{}^t\xi\scirc
\pi_{S,i-2}\scirc{}^t[X] +\delta_{i,3}\,\pi_{f,0},$$
where $\delta_{i,3}=1$ if $i=3$, and $0$ otherwise. Then we have
isomorphisms of Chow motives
$$(X,[X]\scirc\pi_{S,i}\scirc\xi)=(S,\pi_{S,i}),\quad
(X,{}^t\xi\scirc\pi_{S,i-2}\scirc{}^t[X])=(S,\pi_{S,i-2})(-1),$$
using $\xi\scirc[X]=id$ as in (2.4.1--2).
Put $M_{0,j} = (X,\pi_{f,0,j})$.
If $j>r$, we obtain using duality
$$
H^i(M_{0,j})\cong H^{i-2}(C_j,(\iota_j)_*L_j)(-1) = 0
$$
for all $i\ne 3$ in case $\ok=k$, hence
the motive $(X,\pi_{X,3})$ only has cohomology in degree 3.
So, using (1.11) for $j\le r$, we get the Chow-K\"unneth
decomposition for $X$ as desired.

\medskip\noindent
{\bf 2.8.~Proof of Corollary~1.} Using the action of correspondences
on the Chow groups together with (2.4.4), we get
$$\CH_{\alg}^2(X)_{\Q}=\mopl_{-1\le i\le 1}\,
(\pi_{f,i})_*\CH_{\alg}^2(X)_{\Q},$$ and
$$(\pi_{f,-1})_*\CH_{\alg}^2(X)_{\Q}=\CH_{\alg}^2(S)_{\Q},\quad
(\pi_{f,1})_*\CH_{\alg}^2(X)_{\Q}=\CH_{\alg}^1(S)_{\Q},$$ since
$\xi\scirc[X]=id$ as in (2.4.1). We have moreover
$$(\pi_{f,0})_*\CH_{\alg}^2(X)_{\Q}=\mopl_j\,(\tpi_j)_*\CH^1_{\alg}
(D_j)_{\Q} = \mopl_j \CH^1_{\rm alg}(D_j)^{\sigma_j=-1}_{\Q},$$
where the last term is the $(-1)$-eigenspace of
$\CH^1_{\rm alg}(D_j)_{\Q}$ for the action of $\sigma_j$.
So the assertion is reduced to
$$\CH^1_{\alg}(D_j)_{\Q}=J(D_j)(k)_{\Q},$$
where $J(D_j)(k)$ is the abelian group of the $k$-valued points of
the Picard variety of $D_j/k$.
But this is well-known in case $D_j$ has a $k$-valued point,
and the general case is reduced to this case using the action of the
Galois group and the group structure of the Picard variety.
This finishes the proof of Corollary~1.

\medskip\noindent{\bf 2.9.~\bf Relation with Murre's conjectures.}
Let $T(S)\subset\CH^2_{\rm alg}(S)$ be the Albanese kernel, and put
$h^i(S) = (S,\pi_{S,i})$. Assume $S/k$ is absolutely irreducible.
Recall [18] that the rational Chow groups of
the motives $h^i(S)$ are given by the table
$$
\matrix{
& h^0(S) & h^1(S) & h^2(S) & h^3(S) & h^4(S) \cr
\CH^0 & \Q & 0 & 0 & 0 & 0 \cr
\CH^1 & 0 & {\Pic}^0_{S/k}(k)_{\Q} & {\NS}(S)_{\Q} & 0 & 0 \cr
\CH^2 & 0 & 0 & T(S)_{\Q} & {\Alb}_{S/k}(k)_{\Q} & \Q.
}
$$
Assume $r=0$ for simplicity.
Put $M_0 = (X,\pi_{f,0})$, and set $h^i(X) = (X,\pi_{X,i})$. Then
$$\hskip-10pt h(X)\cong h(S)\oplus h(S)(-1)\oplus M_0,$$
and more precisely
$$\eqalign{
&h^0(X)\cong h^0(S),\cr
&h^1(X)\cong h^1(S),\cr
&h^2(X)\cong h^0(S)(-1)\oplus h^2(S),\cr
&h^3(X)\cong h^1(S)(-1)\oplus h^3(S)\oplus M_0,\cr
&h^4(X)\cong h^2(S)(-1)\oplus h^4(S),\cr
&h^5(X)\cong h^3(S)(-1),\cr
&h^6(X)\cong h^4(S)(-1).\cr}$$
Hence the rational Chow groups of the motives $h^i(X)$ are given by the
table
$$
\matrix{
& h^0(X) & h^1(X) & h^2(X) & h^3(X) & h^4(X) & h^5(X) & h^6(X) \cr
\CH^0 & \Q & 0 & 0 & 0 & 0 & 0 & 0 \cr
\CH^1 & 0 & {\Pic}^0_{S/k}(k)_{\Q} & \Q\oplus{\NS}(S)_{\Q} & 0 & 0 & 0 &
0 \cr
\CH^2 & 0 & 0 & T(S)_{\Q} & A_{\Q} & {\NS}(S)_{\Q}\oplus\Q & 0 & 0 \cr
\CH^3 & 0 & 0 & 0 & 0 & T(S)_{\Q} & {\Alb}_{S/k}(k)_{\Q} & \Q
}
$$
with
$$
A_{\Q} = {\Pic}^0_{S/k}(k)_{\Q}\oplus{\Alb}_{S/k}(k)_{\Q}\oplus
\cP_X(k)_{\Q}.
$$
The above table shows that the only correspondences that act nontrivially
on $\CH^j(X)_{\Q}$ are $\pi_{X,j},\ldots,\pi_{X,2j}$.
Hence Murre's conjectures A and B [18] hold for the conic bundle $X$.
This is a refinement of results of del Angel and M\"uller--Stach for
uniruled threefolds [1]. At present, it is not clear whether $X$
satisfies Murre's conjectures C and D.

\medskip\noindent{\bf 2.10.~\bf Remark.}
The decomposition
$$
h(X)\cong h(S) \oplus h(S)(-1) \oplus (\mopl_j \Prym(D_j/C'_j)(-1))
$$
implies that the motive $h(X)$ is finite dimensional (in the sense of
Kimura--O'Sullivan) if $h(S)$ is finite dimensional.

\bigskip\bigskip
\centerline{{\bf References}}

\medskip
{\mfont
\item{[1]}
P.~L.~del Angel and S.~M\"uller-Stach, Motives of uniruled $3$-folds,
Compos. Math. 112 (1998), 1--16.

\item{[2]}
R.~Barlow, Rational equivalence of zero cycles for some more surfaces
with $ p_{g} = 0 $, Inv. Math. 79 (1985), 303--308.

\item{[3]}
A.~Beauville, Vari\'et\'es de Prym et jacobiennes interm\'ediaires,
Ann.\ Sci.\ Ecole Norm.\ Sup.\ (4) 10 (1977), 309--391.

\item{[4]}
A.~Beilinson, Height pairing between algebraic cycles, Lect. Notes
in Math., vol. 1289, Springer, Berlin, 1987, pp. 1--26.

\item{[5]}
A.~Beilinson, J.~Bernstein and P.~Deligne, Faisceaux pervers,
Ast\'erisque, vol. 100, Soc. Math. France, Paris, 1982.

\item{[6]}
M.~Beltrametti, On the Chow group and the intermediate Jacobian of a
conic bundle, Ann. Mat. Pura Appl. (4) 141 (1985), 331--351.

\item{[7]}
S.~Bloch, Lectures on algebraic cycles, Duke University Mathematical
series 4, Durham, 1980.

\item{[8]}
S. Bloch, A. Kas and D. Lieberman, Zero cycles on surfaces with
$ p_{g} = 0, $ Compos. Math. 33 (1976), 135--145.

\item{[9]}
A.~Corti and M.~Hanamura, Motivic decomposition and intersection
Chow groups, I, Duke Math. J. 103 (2000), 459--522.

\item{[10]}
P.~Deligne, Th\'eorie de Hodge II, Publ. Math. IHES,
40 (1971), 5--57.

\item{[11]}
P.~Deligne, Th\'eor\`eme de finitude en cohomologie {\it l}-adique,
in SGA 4 1/2, Lect. Notes in Math., vol. 569, Springer, Berlin, 1977,
233--261.

\item{[12]}
C.~Deninger and J.~P.~Murre, Motivic decomposition of abelian schemes
and the Fourier transform, J. Reine Angew. Math. 422 (1991), 201--219.

\item{[13]}
W.~Fulton, Intersection theory, Springer, Berlin, 1984.

\item{[14]}
B.~B.~Gordon, M.~Hanamura and J.~P.~Murre, Relative Chow-K\"unneth
projectors for modular varieties, J. Reine Angew. Math. 558 (2003),
1--14.

\item{[15]}
B.~B.~Gordon, M.~Hanamura and J.~P.~Murre, Absolute Chow-K\"unneth
projectors for modular varieties, J. Reine Angew. Math. 580 (2005),
139--155.

\item{[16]}
D.~Mumford, Rational equivalence of $0$-cycles on surfaces, J. Math.
Kyoto Univ. 9 (1969), 195--204.

\item{[17]}
J.~P.~Murre, On the motive of an algebraic surface, J. Reine Angew.
Math. 409 (1990), 190--204.

\item{[18]}
J.~P.~Murre, On a conjectural filtration on Chow groups of an
algebraic variety, Indag. Math. 4 (1993), 177--201.

\item{[19]}
J.~Nagel, On the motivic decomposition conjecture for conic bundles,
preprint.

\item{[20]}
M.~Saito, Mixed Hodge modules, Publ.\ RIMS, Kyoto Univ. 26
(1990), 221--333.

\item{[21]}
M.~Saito, Hodge conjecture and mixed motives, I, Proc. Sympos. Pure
Math., 53, Amer. Math. Soc., Providence, RI, 1991, pp.~283--303.

\item{[22]}
M.~Saito, Chow-K\"unneth decomposition for varieties with low
cohomological level, preprint (math.AG/0604254).

\item{[23]}
J.-P.~Serre, Alg\`ebre locale, multiplicit\'es, Lect. Notes in Math.
11, Springer, Berlin, 1975.

\medskip
Jan Nagel

Universit\'e de Lille 1, D\'epartement de Math\'ematiques,
B\^atiment M2

59655 Villeneuve d'Ascq Cedex, France

e-mail: Jan.Nagel@math.univ-lille1.fr

\smallskip
Morihiko Saito

RIMS Kyoto University, Kyoto 606-8502 Japan

e-mail: msaito@kurims.kyoto-u.ac.jp

\medskip
\vers
}
\bye